\documentclass[twoside, 12pt]{article}
\usepackage{amssymb, amsmath, mathrsfs, amsthm}
\usepackage{graphicx}
\usepackage{color}
\usepackage[top=2cm, bottom=2cm, left=4cm, right=4cm]{geometry}
\usepackage{float}

\newcommand{\bd}{\begin{description}}
\newcommand{\ed}{\end{description}}
\newcommand{\bi}{\begin{itemize}}
\newcommand{\ei}{\end{itemize}}
\newcommand{\be}{\begin{enumerate}}
\newcommand{\ee}{\end{enumerate}}
\newcommand{\beq}{\begin{equation}}
\newcommand{\eeq}{\end{equation}}
\newcommand{\beqs}{\begin{eqnarray*}}
\newcommand{\eeqs}{\end{eqnarray*}}

\definecolor{DarkGreen}{rgb}{0.2, 0.6, 0.3}

\newtheorem{theorem}{Theorem}
\newtheorem{conjecture}{Conjecture}

\newtheorem{case}{Case}

\newtheorem{claim}{Claim}

\begin{document}
\title{\textbf{A note on a new result related to Chartrand, Kaigars and Lick's theorem }}
\author{Zhong Huang\footnote{College of Information and Mathematics,
Yangtze University,
Jingzhou, China. {\tt hz@yangtzeu.edu.cn}},\ \ Meng Ji\footnote{Corresponding author: College of Mathematical Science, Tianjin Normal University,
Tianjin, China.} \footnote{Institute of Mathematics and Interdisciplinary Sciences, Tianjin Normal University, Tianjin, China. {\tt mji@tjnu.edu.cn}}}

\date{}
\maketitle

\begin{abstract}
In this note, we prove a theorem covering Chartrand, Kaigars, and Lick's theorem in [Proc. Amer. Math. Soc. 32 (1972), 63–68]. As an application, we give a simpler proof of theorem proved by Mader [J. Graph Theory 65 (2010), 61--69. (Theorem 1)].
\\[2mm]
\end{abstract}

\section{Introduction}

In this paper, all the graphs are finite, undirected and simple. For
graph-theoretical terminology and notations not defined here, we
follow \cite{BM}.
Let $G=(V(G),E(G))$ be a graph.
For any vertex $x\in V(G)$, we denote by $N_{G}(x)$ the set of neighbors of $x$ in $G$. 
We denote by $P_n$ the path of order $n$. 
A graph $G$ is said to be \emph{$k$-connected} if the removal of fewer than $k$ vertices from $G$ neither disconnects it nor reduces it to the trivial graph consisting of a single vertex. 
The maximum value of $k$ for which a graph $G$ is $k$-connected is called its \emph{connectivity} and is denoted by $\kappa(G)$.
A $k$-connected graph is \emph{critically $k$-connected}, if deleting any vertex, the connectivity becomes less than $k$.
The \emph{minimum degree} of $G$ is denoted by $\delta(G)$.
A \emph{separating set} of a connected graph $G$ is a subset $S$ of the vertex set $V(G)$ such that $G-S$ has more than one connected component.
And a \emph{minimum separating set} is a separating set with smallest size in a given connected graph $G$, and the size of a minimum separating set is $\kappa(G)$.
For a minimum separating set $S$ of a graph $G$, a union of connected components $F$ of $G-S$, with
$F\ne G-S$, is a \emph{fragment} $F$ to $S$, and the
\emph{complementary fragment} $G-(S\cup V(F))$ is denoted by $\bar{F}$.
If a fragment of $G$ does not contain any other fragments of $G$ to certain minimum separating sets, then
it is an \emph{end} of $G$.
Clearly, every graph contains an end except for complete graphs.

In 1972, Chartrand, Kaigars and Lick \cite{CKL} proved that \emph{every critically
$k$-connected graph contains a vertex of degree less than $\lfloor\frac{3k}{2}\rfloor$.}

\begin{theorem}[\upshape\bf{Chartrand, Kaigars and Lick} \cite{CKL}] \label{Chartrand-Kaigars-Lick}
Every $k$-connected graph $G$ with $\delta (G)\geq
\lfloor\frac{3k}{2}\rfloor$ has a vertex $x$ with $\kappa(G-x)\geq k$.
\end{theorem}

Here we give a stronger theorem as follows.

\begin{theorem}\label{them1}
For every positive integer $k$, let $G$ be a graph with $\kappa(G)=k$ and $\delta(G)\geq\lfloor\frac{3k}{2}\rfloor$, let fragment $F$ to a minimum separating set $S$ be an end of $G$. Then for each minimum separating set $S^\prime$ of $G$, we have that $F\cap S^\prime=\emptyset$.
\end{theorem}

Let $F$ be an end of graph $G$.
Then the removal of any vertex $v\in F$ will not decrease the connectivity $\kappa(G)$ by Theorem \ref{them1}.
We have that Theorem \ref{Chartrand-Kaigars-Lick} is a immediate corollary of our theorem.

In 2008, Fujita and Kawarabayashi \cite{SK} showed the following result.

\begin{theorem}[{\bf Fujita and Kawarabayashi} \cite{SK}]
Every $k$-connected graph $G$ with
$\delta (G)\geq \lfloor\frac{3k}{2}\rfloor+2$ has an edge $xy$ such
that $G-\{x,y\}$ remains $k$-connected.
\end{theorem}

Furthermore, Fujita and Kawarabayashi posed the following conjecture.

\begin{conjecture}[{\bf Fujita and Kawarabayashi} \cite{SK}]\label{fujita}
For all positive integers $k$, $m$, there is a (least) non-negative
integer $f_{k}(m)$ such that every $k$-connected graph $G$ with
$\delta(G)\geq \lfloor\frac{3k}{2}\rfloor+ f_{k}(m)-1$ contains a
connected subgraph $W$ of order $m$ such that $G-V(W)$ is still
$k$-connected.
\end{conjecture}

In 2010, Mader \cite{mader1} proved that Conjecture \ref{fujita} holds by providing $f_{k}(m)=m$ and $W=P_m$.

\begin{theorem}[{\bf Mader} \cite{mader1}]\label{made}
Every $k$-connected graph $G$ with $\delta(G)\geq \lfloor\frac{3k}{2}\rfloor+ m-1$ for positive integers $k,m$
contains a path $P$ of $|P|=m$ such that $G-V(P)$ is still $k$-connected.
\end{theorem}

As another application of Theorem \ref{them1}, we give a simpler proof of Theorem \ref{made}.

\section{Proof of Theorem \ref{them1}}

Let fragment $F$ to a minimum separating set $S$ be an end of $G$.
Without loss of generality, we can assume that $F^\prime$ is the fragment to a minimum separating set $S^\prime$ of $G$ such that $|V(F^\prime)\cap V(F)|\ge |V(\bar{F^\prime})\cap V(F)|$.
The vertex set of $G$ is partitioned into three disjoint sets $V(F)$, $S$ and $V(\bar{F})$.
The vertex set $V(F)$ is partitioned into three disjoint sets $V(F^\prime)\cap V(F)$, $S^\prime\cap V(F)$ and $V(\bar{F^\prime})\cap V(F)$, denoted by $A_1$, $A_2$ and $A_3$, respectively.
The vertex set $S$ is partitioned into three disjoint sets $V(F^\prime)\cap S$, $S^\prime\cap S$ and $V(\bar{F^\prime})\cap S$, denoted by $B_1$, $B_2$ and $B_3$, respectively.
The vertex set $V(\bar{F})$ is partitioned into three disjoint sets $V(F^\prime)\cap V(\bar{F})$, $S^\prime\cap V(\bar{F})$ and $V(\bar{F^\prime})\cap V(\bar{F})$, denoted by $C_1$, $C_2$ and $C_3$, respectively.
Clearly, we have that $V(F^\prime)= A_1\cup B_1\cup C_1$, $S^\prime= A_2\cup B_2\cup C_2$, $V(\bar{F^\prime})= A_3\cup B_3\cup C_3$ and $|A_1|\geq |A_3|$ by assumption.
Next we prove the theorem by contradiction.
Suppose that $F\cap S^\prime\ne \emptyset$.
\setcounter{case}{0}
\begin{case}
 $|A_1|=0$.
\end{case}
\noindent Since $|A_1|\geq |A_3|$, both $A_1$ and $A_3$ are empty sets, furthermore, we have $|A_2|> 0$. So there is a vertex $v\in A_2$.
Since $N_G(v)\subset (A_2\cup S)\setminus \{v\}$, $|S|=k$ and $|N_G(v)|\ge \lfloor\frac{3k}{2}\rfloor$, we have that $|A_2|\ge \lfloor\frac{k}{2}\rfloor +1$.
Similarly, if $|C_1|=0$, then $|B_1|\ge \lfloor\frac{k}{2}\rfloor +1$.
And if $|C_3|=0$, then $|B_3|\ge \lfloor\frac{k}{2}\rfloor +1$.
Suppose that $|C_1|=|C_3|=0$, then $|S|\ge |B_1|+|B_3|\ge 2\lfloor\frac{k}{2}\rfloor +2 >k$, a contradiction. 
Suppose that $|C_1|=0$ and $|C_3|>0$, then $C_2\cup B_2 \cup B_3$ is a separating set such that $|B_3\cup B_2 \cup C_2|\le |C_2|+2|B_2|+|B_3|=|S|-|B_1|+|S^\prime|-|A_2|\le 2k-2\lfloor\frac{k}{2}\rfloor-2<k$, a contradiction.
Suppose that $|C_1|> 0$ and $|C_3|=0$, then $C_2\cup B_2 \cup B_1$ is a separating set such that $|B_1\cup B_2 \cup C_2|<k$ by similarly argument, a contradiction.
Suppose that $|C_1|> 0$ and $|C_3|> 0$, then both $B_3\cup B_2 \cup C_2$ and $B_1\cup B_2 \cup C_2$ are separating sets in $G$.
Since $|B_3\cup B_2 \cup C_2|+|B_1\cup B_2 \cup C_2|\le |B_1|+3|B_2|+|B_3|++2|C_2|=|S|+2(|S^\prime|-|A_2|)\le 3k-2\lfloor\frac{k}{2}\rfloor-2<2k $, we have that either $|B_3\cup B_2 \cup C_2|<k$ or $|B_1\cup B_2\cup C_2|<k$ by the pigeonhole principle, a contradiction.

\begin{case}
 $|A_1|>0$ and $|A_3|=0$.
\end{case}
\noindent It is easy to know that $B_1\cup B_2\cup A_2$ is a separating set.
We have that $|B_1\cup B_2\cup A_2|>k$, otherwise, $F$ is not an end to a minimum separating set $S$ of $G$, a contradiction.
Suppose that $|C_3|=0$, then both $A_3$ and $C_3$ are empty sets.
By similar argument in Case 1, we deduce a contradiction.
Suppose that $|C_3|> 0$, then $C_2\cup B_2\cup B_3$ is a separating set such that $| B_3\cup B_2\cup C_2|=|B_3|+|B_2|+|C_2|=(|S^\prime|-|A_2|)+(|S|-|B_1|)-|B_2|=2k-|B_1\cup B_2\cup A_2|<k$, a contradiction.

\begin{case}
 $|A_1|>0$ and $|A_3|>0$.
\end{case}
\noindent Both $B_1\cup B_2\cup A_2$ and $B_3\cup B_2\cup A_2$ are separating sets.
We have that $|B_1\cup B_2\cup A_2|>k$ and $|B_3\cup B_2\cup A_2|>k$, otherwise $F$ is not an end to a minimum separating set $S$ of $G$, a contradiction.
Suppose that $|C_1|> 0$, then $B_1\cup B_2\cup C_2$ is a separating set such that $|B_1\cup B_2\cup C_2|=|B_1|+|B_2|+|C_2|=(|S^\prime|-|A_2|)+(|S|-|B_3|)-|B_2|=2k-|B_3\cup B_2\cup A_2|<k$, a contradiction.
Suppose that $|C_3|> 0$, similarly we deduce a contradiction.
Suppose that $|C_1|=|C_3|=0$, then both $C_1$ and $C_3$ are empty sets.
By similar argument in Case 1, we deduce a contradiction.
This completes the proof.
\qed \\[0.2mm]

\section{Simpler proof of Theorem \ref{made}}
For a graph $G$ and any vertex $v\in G$, it follows that $\kappa(G)\ge \kappa(G-v)$.
\setcounter{claim}{0}
\begin{claim}\label{claim1}
There exists a set of vertices $\{v_{1},v_{2},\cdots ,v_{m}\}$ in $G$ ($G_0=G$) satisfying the following properties for $0\le i\le m-1$. 
$i)$ $v_{i+1}$ is adjacent to $v_{i}$ in $G$;
$ii)$ the subgraph $G_{i+1}=G_i-v_{i+1}$ is $k$-connected;
$iii)$ if $\kappa(G_{i})\ge\kappa(G_{i+1})>k$, then $F_i=G_i$ and $F_{i+1}=G_{i+1}$;
 if $\kappa(G_{i})>\kappa(G_{i+1})=k$, then $F_i=G_i$ and there exists an end $F_{i+1}$ in $G_{i+1}$;
 if $\kappa(G_{i})=\kappa(G_{i+1})=k$, then $F_i$ is an end of $G_i$ and there exists an end $F_{i+1}$ in $G_{i+1}$;
$iv)$ $F_{i+1}\cap N_G(v_{i+1})\ne \emptyset$.
\end{claim}

\begin{proof}
We prove by induction.
At the beginning for $i=0$.
Let $F_0=G_0$, we pick an arbitrary vertex $v_1\in F_0$.
Since $v_0$ does not exist, we have that $i)$ trivially holds by definition for supplementary.
By Theorem \ref{them1}, we have that $G_1=G_0-v_1$ is $k$-connected, $ii)$ holds.
Suppose that $\kappa(G_{0})\ge\kappa(G_{1})>k$, then $F_0=G$, $F_{1}=G_{1}$ and $F_{1}\cap N_G(v_{1})\ne \emptyset$.
Suppose that $\kappa(G_{0})>\kappa(G_{1})=k$, then $F_0=G$ and there exists an end $F_{1}$ in $G_{1}$ by definition. 
We have that $F_{1}\cap N_G(v_{1})\ne \emptyset$, otherwise $F_1$ is also an end of $G_0$ such that $\kappa(G_{0})=k$, a contradiction.
Suppose that $\kappa(G_{0})=\kappa(G_{1})=k$, then $F_{0}$ is an end in $G_{0}$.
If $F_0-v_1$ is an end of $G_1$, then we pick $F_1=F_0-v_1$.
If $F_0-v_1$ is not an end of $G_1$, then $F_0-v_1$ is a fragment of $G_1$ and by the definition of the end we pick an end $F_1$ in $G_1$ such that $F_1\subset F_0-v_1\subset F_0$.
We have that $F_{1}\cap N_G(v_{1})\ne \emptyset$, otherwise $F_1\subset F_0$ and $F_1$ is also an end of $G_0$, a contradiction.
$iii)$ and $iv)$ hold.
Consequently, Claim 1 holds for $i=0$. 
Assume that Claim 1 holds for $i=t$ with $0\le t\le m-2$.
It follows that $F_{t+1}\cap N_G(v_{t+1})\ne \emptyset$ by assumption.
We pick a vertex $v_{t+2}\in F_{t+1}\cap N_G(v_{t+1})$.
Clearly, $v_{i+2}$ is adjacent to $v_{i+1}$, $i)$ holds for $i=t+1$.
Since $\delta(G_{t+1})\ge \delta(G)-(t+1)\ge \lfloor\frac{3k}{2}\rfloor -t+m-2\ge \lfloor\frac{3k}{2}\rfloor$,
then we have that $G_{t+2}$ is $k$-connected by Theorem \ref{them1}, $ii)$ holds for $i=t+1$.
Suppose that $\kappa(G_{t+1})\ge\kappa(G_{t+2})>k$, then $F_{t+2}=G_{t+2}$.
Since $N_G(v_{t+2})> m\ge t+2$, then we have that $F_{t+2}\cap N_G(v_{t+2})\ne \emptyset$.
Suppose that $\kappa(G_{t+1})>\kappa(G_{t+2})=k$, then $F_{t+1}=G_{t+1}$ and there exists an end $F_{t+2}$ in $G_{t+2}$ by definition.
We have that $F_{t+2}\cap N_G(v_{t+2})\ne \emptyset$, otherwise $F_{t+2}$ is also an end of $G_{t+1}$ such that $\kappa(G_{t+1})=k$, a contradiction.
Suppose that $\kappa(G_{t+1})=\kappa(G_{t+2})=k$, then $F_{t+1}$ is an end in $G_{t+1}$.
If $F_{t+1}-v_{t+2}$ is an end of $G_{t+2}$, then we pick $F_{t+2}=F_{t+1}-v_{t+2}$.
If $F_{t+1}-v_{t+2}$ is not an end of $G_{t+2}$, then $F_{t+1}-v_{t+2}$ is a fragment of $G_{t+2}$ and we pick an end $F_{t+2}$ in $G_{t+2}$ such that $F_{t+2}\subset F_{t+1}-v_{t+2}\subset F_{t+1}$.
We have that $F_{t+2}\cap N_G(v_{t+2})\ne \emptyset$, otherwise $F_{t+2}\subset F_{t+1}$ and $F_{t+2}$ is also an end of $G_{t+1}$, a contradiction.
It follows that $iii)$ and $iv)$ hold.
Then Claim 1 holds for $i=t+1$.
\end{proof}

Let $i=m-1$, then there exists a path $P=v_1v_2\cdots v_m$ such that $G-V(P)$ is $k$-connected by Claim 1.
This completes the proof.

\qed \\[0.2mm]

\section{Concluding remarks}
In 2009, Mader \cite{mader1} posed a conjecture which is open for $k\geq 4$. Theorem \ref{them1} provides a new method to deal with this conjecture by constructing connectivity keeping subgraphs.
\begin{conjecture}[{\bf Mader} \cite{mader1}]\label{conj}
For every positive integer $k$ and finite tree $T$ of order $m$,
every $k$-connected finite graph $G$ with minimum degree
$\delta(G)\geq\lfloor\frac{3k}{2}\rfloor+m-1$ contains a subgraph
$T'\cong T$ such that $G-V(T')$ remains $k$-connected.
\end{conjecture}

\noindent {\bf Declaration of Interest Statement}: The authors declare no conflict of interest.

\noindent {\bf Acknowledgement}: We would like
to thank the anonymous reviewers for the valuable suggestions and comments for the previous version.


\begin{thebibliography}{1}

\bibitem{BM}
J. A. Bondy and U. S. R. Murty, \emph{Graph Theory}, Graduate Texts
in Mathematics 244, Springer, Berlin, 2008.

\bibitem{CKL}
G. Chartrand, A. Kaigars, and D.R. Lick, Critically $n$-connected
graphs, \emph{Proc. Amer. Math. Soc.} 32 (1972), 63--68.

\bibitem{SK}
S. Fujita and K. Kawarabayashi, Connectivity keeping edges in graphs
with large minimum degree, \emph{J. Combin. Theory, Ser. B} 98
(2008), 805--811.

\bibitem{mader1}
W. Mader, Connectivity keeping paths in $k$-connected graphs,
\emph{J. Graph Theory} 65 (2010), 61--69.

\end{thebibliography}
\end{document}